\documentclass[10pt]{amsart}
\usepackage{amsfonts}
\usepackage{ifthen}
\usepackage{amsthm}
\usepackage{amsmath}
\usepackage{graphicx}
\usepackage{amscd,amssymb,amsthm}

\newcounter{minutes}
\divide\time by 60
\newcounter{hours}
\multiply\time by 60 \addtocounter{minutes}{-\time}

\newcommand{\real}{\operatorname{Re}}

\newtheorem{theorem}{Theorem}
\newtheorem{lemma}{Lemma}


\title[]{Starlikeness of Bessel functions and their derivatives}

\author[\'A. Baricz]{\'Arp\'ad Baricz$^{\bigstar}$}
\address{Department of Economics, Babe\c{s}-Bolyai University, Cluj-Napoca 400591, Romania}
\address{Institute of Applied Mathematics, John von Neumann Faculty of Informatics, \'Obuda University, 1034 Budapest, Hungary}
\email{bariczocsi@yahoo.com}

\author[{M. \c{C}a\u{g}lar}]{Murat \c{C}a\u{g}lar}
\address{Department of Mathematics, Faculty of Science and Letters, Kafkas
University, Kars 36100, Turkey.} \email{mcaglar25@gmail.com}

\author[E. Deniz]{Erhan Deniz}
\address{Department of Mathematics, Faculty of Science and Letters, Kafkas
University, Kars 36100, Turkey.} \email{edeniz36@gmail.com}

\thanks{$^{\bigstar}$The research of \'A. Baricz was supported by a research grant of the Romanian National Authority for Scientific Research, CNCS-UEFISCDI, project number PN-II-RU-TE-2012-3-0190. The work of \'A. Baricz was completed during his visit in September 2014 to Department of Mathematics of Kafkas University in Kars, to which this author is grateful for hospitality.}

\begin{document}

\def\thefootnote{}
\footnotetext{ \texttt{File:~\jobname .tex,
          printed: \number\year-\number\month-0\number\day,
          \thehours.\ifnum\theminutes<10{0}\fi\theminutes}
} \makeatletter\def\thefootnote{\@arabic\c@footnote}\makeatother

\maketitle

\begin{center}
{Dedicated to Aysima, Bor\'oka, Eren and Kopp\'any}
\end{center}

\begin{abstract}
In this paper necessary and sufficient conditions are deduced for the starlikeness of Bessel functions of the first kind and their derivatives of the second and third order by using a result of Shah and Trimble about transcendental entire functions with univalent derivatives and some Mittag-Leffler expansions for the derivatives of Bessel functions of the first kind, as well as some results on the zeros of these functions.
\end{abstract}

\section{\bf Introduction and the Main Results}

Geometric properties of Bessel functions of the first kind $J_{\nu}$, like univalence, starlikeness, spirallikeness and convexity were studied in the sixties by Brown \cite{brown, brown2,brown3}, and also by Kreyszig and Todd \cite{todd}. Other geometric properties of Bessel functions of the first kind were studied later in the papers \cite{mathematica,publ,lecture,bsk,samy,basz,szasz,szasz2}. Very recently, in \cite{barsza} the close-to-convexity of the derivatives of Bessel functions was considered. Motivated by the above results, in this paper we make a contribution to the subject by obtaining some necessary and sufficient conditions for the starlikeness of Bessel functions of the first kind and their derivatives of the second and third order by using a result of Shah and Trimble \cite[Theorem 2]{st} about transcendental entire functions with univalent derivatives and some Mittag-Leffler expansions for the derivatives of Bessel functions of the first kind, as well as some results on the zeros of these functions.

Our first set of sharp results are about the starlikeness of order $\alpha$ of two normalized Bessel functions of the first kind. We note that these results naturally complement the main results of \cite{bsk,basz,szasz}.

\begin{theorem}\label{th1}
The function
$$z\mapsto f_{\nu}(z)=\left(2^{\nu}\Gamma(\nu+1)J_{\nu}(z)\right)^{\frac{1}{\nu}}$$
is starlike of order $\alpha\in[0,1)$ in $\mathbb{D}$ if and only if $\nu>\nu_1(\alpha),$ where $\nu_1(\alpha)$ is the unique root of the equation
$(1-\alpha)\nu J_{\nu}(1)=J_{\nu+1}(1),$ situated in $(0,\infty).$ In particular, $f_{\nu}$ is starlike in $\mathbb{D}$ if and only if $\nu>\nu_1(0),$ where $\nu_1(0)\simeq0.3908\dots$ is the unique root of the transcendental equation
$\nu J_{\nu}(1)=J_{\nu+1}(1).$
\end{theorem}

\begin{theorem}\label{th2}
The function
$$z\mapsto g_{\nu}(z)=2^{\nu}\Gamma(\nu+1)z^{1-\nu}J_{\nu}(z)$$
is starlike of order $\alpha\in[0,1)$ in $\mathbb{D}$ if and only if $\nu>\nu_2(\alpha),$ where $\nu_2(\alpha)$ is the unique root of the equation
$(1-\alpha)J_{\nu}(1)=J_{\nu+1}(1),$ situated in $(\tilde{\nu},\infty),$ where $\tilde{\nu}\simeq-0.7745\dots$ is the unique root of $j_{\nu,1}=1$ and $j_{\nu,1}$ is the first positive zero of $J_{\nu}.$ In particular, the function $g_{\nu}$ is starlike in $\mathbb{D}$ if and only if $\nu>\nu_2(0),$ where $\nu_2(0)\simeq-0.3397\dots$ is the unique root of the transcendental equation
$J_{\nu}(1)=J_{\nu+1}(1).$
\end{theorem}

We note that very recently Antonino and Miller \cite[Example 3]{antonino} as an application of the third-order differential subordinations proved that the function $z\mapsto \int\limits_0^zJ_0(t)dt$ is convex (and hence univalent) in $\mathbb{D}.$ If we consider the function $w_{\nu}:\mathbb{D}\to\mathbb{C},$ defined by
$$w_{\nu}(z)=2^{\nu}\Gamma(\nu+1)\int_0^zt^{-\nu}J_{\nu}(t)dt,$$
then in view of the relation
$$1+\frac{zw_{\nu}''(z)}{w_{\nu}'(z)}=\frac{zg_{\nu}'(z)}{g_{\nu}(z)}$$
and the analytic characterizations of starlike and convex functions, Theorem \ref{th2} can be rewritten as follows: the function $w_{\nu}$ is convex of order $\alpha\in[0,1)$ in $\mathbb{D}$ if and only if $\nu>\nu_2(\alpha),$ and in particular, the function $w_{\nu}$ is convex (and hence univalent) in $\mathbb{D}$ if and only if $\nu>\nu_2(0).$ This generalizes the result of Antonino and Miller \cite[Example 3]{antonino} on $w_0$ and shows actually that if $\nu<\nu_2(0),$ then the above convexity property is no more true.

The next set of sharp main results are based on a result of Shah and Trimble \cite[Theorem 2]{st}, see Lemma \ref{lem} in the next section, and these results are natural companions of the main results in \cite{bdy,barsza}. We note that it would be interesting to see a common generalization of the next three theorems. Following the proof of these theorems it is clear that the monotonicity of the zeros (with respect to the order) of the derivative (of arbitrary order greater than three) of Bessel functions of the first kind would be enough together with Lemma \ref{lem}.

\begin{theorem}\label{thA}
The function $$z\mapsto 2^{\nu}\Gamma(\nu)z^{\frac{3}{2}-\frac{\nu}{2}}J_{\nu}'(\sqrt{z})$$ is starlike and all of its derivatives are close-to-convex (and hence univalent) in $\mathbb{D}$ if and only if $\nu\geq\dot{\nu},$ where $\dot{\nu}\simeq0.7022\dots$ is the unique root on $(0,\infty)$ of the transcendent equation  $$(2\nu-1)J_{\nu}(1)+(\nu-2)J_{\nu+1}(1)=0.$$
\end{theorem}

\begin{theorem}\label{thB}
The function $$z\mapsto 2^{\nu}\Gamma(\nu-1)z^{2-\frac{\nu}{2}}J_{\nu}''(\sqrt{z})$$ is starlike and all of its derivatives are close-to-convex (and hence univalent) in $\mathbb{D}$ if and only if $\nu\geq\nu^{\ast},$ where $\nu^{\ast}\simeq1.9052\dots$ is the unique root on $(1,\infty)$ of the transcendent equation  $$(2\nu^2-2\nu-3)J_{\nu}(1)=(\nu^2+\nu-3)J_{\nu+1}(1).$$
\end{theorem}

\begin{theorem}\label{thC}
The function $$z\mapsto 2^{\nu}\Gamma(\nu-2)z^{\frac{5}{2}-\frac{\nu}{2}}J_{\nu}'''(\sqrt{z})$$ is starlike and all of its derivatives are close-to-convex (and hence univalent) in $\mathbb{D}$ if and only if $\nu\geq\nu^{\star},$ where $\nu^{\star}\simeq3.077\dots$ is the unique root on $(2,\infty)$ of the transcendent equation  $$(2\nu^3-7\nu^2+3)J_{\nu}(1)+(\nu^3+\nu^2+\nu-1)J_{\nu+1}(1)=0.$$
\end{theorem}

The last main result of this paper is a common generalization of Theorems \ref{thA} and \ref{thB}.

\begin{theorem}\label{thD}
Let $a,b,c\in\mathbb{R}$ such that $c=0$ and $b\neq a$ or $c>0$ and $b>a.$ Moreover, suppose that $\nu\geq \overline{\nu},$ where $\overline{\nu}=\max\{0,\nu_0\}$ and $\nu_0$ is the largest root of the quadratic $Q(\nu)=a\nu(\nu-1)+b\nu+c.$ Assume also that the following inequalities are valid
\begin{equation}\label{abc}
Q(\nu)+4a\nu +2a+2b>0,\ \ \ (4\nu+3)Q(\nu)>4a\nu+2a+2b.
\end{equation}
Then the function $$z\mapsto 2^{\nu}Q^{-1}({\nu})\Gamma(\nu+1)z^{1-\frac{\nu}{2}}\left(azJ_{\nu}''(\sqrt{z})+b\sqrt{z}J_{\nu}'(\sqrt{z})+cJ_{\nu}(\sqrt{z})\right)$$ is starlike and all of its derivatives are close-to-convex (and hence univalent) in $\mathbb{D}$ if and only if $\nu\geq\nu^{\circ},$ where $\nu^{\circ}$ is the unique root on $(\overline{\nu},\infty)$ of the transcendent equation $$(2a\nu^2-2a\nu+2b\nu-3a-b+2c)J_{\nu}(1)=(a\nu^2+a\nu-b\nu-3a+2b+c)J_{\nu+1}(1).$$
\end{theorem}

It is worth to mention that when $b=c=0$ and $a=1,$ then Theorem \ref{thD} reduces to Theorem \ref{thB}. In this case $\overline{\nu}=1,$ $\nu^{\circ}$ becomes $\nu^{\star}$ and the inequalities \eqref{abc} become $\nu^2+3\nu+2>0,$ and $4\nu^3-\nu^2-7\nu-2>0.$ These inequalities give $\nu>-1$ and $\nu>1.5687{\dots},$ which are certainly satisfied for $\nu>\nu^{\star}.$

Similarly, we note that when $a=c=0$ and $b=1,$ then Theorem \ref{thD} reduces to Theorem \ref{thA}. In this case $\overline{\nu}=0,$ $\nu^{\circ}$ becomes $\dot{\nu}$ and the inequalities \eqref{abc} become $\nu+2>0,$ and $4\nu^2+3\nu-2>0.$ These inequalities give $\nu>-2$ and $\nu>0.4253{\dots},$ which are certainly satisfied for $\nu>\dot{\nu}.$

\section{\bf Proofs of the main results}

In this section our aim is to present the proof of the main results of this paper. The proofs of Theorems \ref{th1} and \ref{th2} are mainly based on the Mittag-Leffler expansions and some inequalities from the proof of the main result from \cite{bsk}.

\begin{proof}[\bf Proof of Theorem \ref{th1}]
Let us denote by $j_{\nu,n}$ the $n$th positive zero of the function $J_{\nu}.$ From the proof of \cite[Theorem 1]{bsk} we know that for $\nu>0$ and $r=|z|<j_{\nu,1}$ we have that
$$\real \left(\frac{zf_{\nu}'(z)}{f_{\nu}(z)}\right)\geq \frac{rf_{\nu}'(r)}{f_{\nu}(r)}=1-\frac{1}{\nu}\sum_{n\geq 1}\frac{2r^2}{j_{\nu,n}^2-r^2}.$$
Since $j_{\nu,1}>j_{0,1}\simeq2.4048{\dots}>1$ when $\nu>0,$ the above inequality is clearly valid when $|z|<1.$ On the other hand, the function $r\mapsto rf_{\nu}'(r)/f_{\nu}(r)$ is clearly decreasing on $(0,1)\subset(0,j_{\nu,1}),$ and consequently for all $z\in\mathbb{D}$ and $\nu>0$ we have
$$\real \left(\frac{zf_{\nu}'(z)}{f_{\nu}(z)}\right)\geq \frac{rf_{\nu}'(r)}{f_{\nu}(r)}=1-\frac{1}{\nu}\sum_{n\geq 1}\frac{2r^2}{j_{\nu,n}^2-r^2}\geq1-\frac{1}{\nu}\sum_{n\geq 1}\frac{2}{j_{\nu,n}^2-1}=\frac{f_{\nu}'(1)}{f_{\nu}(1)}.$$
Since the function $\nu\mapsto j_{\nu,n}$ is increasing on $(0,\infty)$ for $n\in\mathbb{N}$ fixed (see \cite[p. 236]{nist}), it follows that the function $\nu\mapsto f_{\nu}'(1)/f_{\nu}(1)$ is increasing on $(0,\infty),$ and thus $f_{\nu}'(1)/f_{\nu}(1)>\alpha$ if and only if $\nu>\nu_1(\alpha),$ where $\nu_1(\alpha)$ is the unique root of the equation
$$f_{\nu}'(1)=\alpha f_{\nu}(1)\ \ \Longleftrightarrow \ \ \nu\alpha J_{\nu}(1)=J_{\nu}'(1) \ \ \Longleftrightarrow\ \ (1-\alpha)\nu J_{\nu}(1)=J_{\nu+1}(1).$$ Here we used that
$$\frac{zf_{\nu}'(z)}{f_{\nu}(z)}=\frac{1}{\nu}\frac{zJ_{\nu}'(z)}{J_{\nu}(z)}=1-\frac{J_{\nu+1}(z)}{\nu J_{\nu}(z)}.$$
Taking into account the fact that all of the above inequalities are sharp it follows that indeed the function $f_{\nu}$ is starlike of order $\alpha\in[0,1)$ in $\mathbb{D}$ if and only if $\nu>\nu_1(\alpha).$
\end{proof}

\begin{proof}[\bf Proof of Theorem \ref{th2}]
From the proof of \cite[Theorem 1]{bsk} we know that for $\nu>-1$ and $r=|z|<j_{\nu,1}$ we have that
$$\real \left(\frac{zg_{\nu}'(z)}{g_{\nu}(z)}\right)\geq \frac{rg_{\nu}'(r)}{g_{\nu}(r)}=1-\sum_{n\geq 1}\frac{2r^2}{j_{\nu,n}^2-r^2}.$$
Since $\nu\mapsto j_{\nu,1}$ is increasing on $(-1,\infty),$ it follows that $j_{\nu,1}>1$ when $\nu>\tilde{\nu},$ and thus in this case the above inequality is clearly valid when $|z|<1.$ On the other hand, the function $r\mapsto rg_{\nu}'(r)/g_{\nu}(r)$ is clearly decreasing on $(0,1)\subset(0,j_{\nu,1}),$ and consequently for all $z\in\mathbb{D}$ and $\nu>\tilde{\nu}$ we have
$$\real \left(\frac{zg_{\nu}'(z)}{g_{\nu}(z)}\right)\geq \frac{rg_{\nu}'(r)}{g_{\nu}(r)}=1-\sum_{n\geq 1}\frac{2r^2}{j_{\nu,n}^2-r^2}\geq1-\sum_{n\geq 1}\frac{2}{j_{\nu,n}^2-1}=\frac{g_{\nu}'(1)}{g_{\nu}(1)}.$$
Since the function $\nu\mapsto j_{\nu,n}$ is increasing on $(-1,\infty)$ for $n\in\mathbb{N}$ fixed (see \cite[p. 236]{nist}), it follows that the function $\nu\mapsto g_{\nu}'(1)/g_{\nu}(1)$ is increasing on $(\tilde{\nu},\infty),$ and thus $g_{\nu}'(1)/g_{\nu}(1)>\alpha$ if and only if $\nu>\nu_2(\alpha),$ where $\nu_2(\alpha)$ is the unique root of the equation
$$g_{\nu}'(1)=\alpha g_{\nu}(1)\ \ \Longleftrightarrow \ \ (1-\nu-\alpha)J_{\nu}(1)+J_{\nu}'(1)=0 \ \ \Longleftrightarrow\ \ (1-\alpha)J_{\nu}(1)=J_{\nu+1}(1).$$ Here we used that
$$\frac{zg_{\nu}'(z)}{g_{\nu}(z)}=1-\nu+\frac{zJ_{\nu}'(z)}{J_{\nu}(z)}=1-\frac{zJ_{\nu+1}(z)}{J_{\nu}(z)}.$$
Taking into account the fact that all of the above inequalities are sharp it follows that indeed the function $g_{\nu}$ is starlike of order $\alpha\in[0,1)$ in $\mathbb{D}$ if and only if $\nu>\nu_2(\alpha).$
\end{proof}

Now, for the proof of the remaining theorems we will use the following result of Shah and Trimble \cite[Theorem 2]{st} about transcendental entire functions with univalent derivatives, which was the key tool in the proof of the main results of \cite{bdy,basz}.

\begin{lemma}\label{lem}
Let $\mathbb{D}=\{z\in\mathbb{C}:|z|<1\}$ be the open unit disk and $f:\mathbb{D}\to\mathbb{C}$ be a transcendental entire function of the form
$$f(z)=z\prod_{n\geq 1}\left(1-\frac{z}{z_n}\right),$$
where all $z_n$ have the same argument and satisfy $|z_n|>1.$ If $f$ is univalent in $\mathbb{D},$ then
$$\sum_{n\geq1}\frac{1}{|z_n|-1}\leq 1.$$
In fact the above inequality holds if and only if $f$ is starlike in $\mathbb{D}$ and all of its derivatives are
close-to-convex there.
\end{lemma}

As we can see below the structures of the next proofs are very similar and all of them use the monotonicity of the zeros with respect to the order of the derivatives of Bessel functions of the first kind.

\begin{proof}[\bf Proof of Theorem \ref{thA}]
Let us denote by $j_{\nu,n}'$ the $n$th positive zero of the function $J_{\nu}'.$ By using the infinite product representation \cite[p. 340]{skelton}
$$J_{\nu}'(z)=\frac{\left(\frac{z}{2}\right)^{\nu-1}}{2\Gamma(\nu)}\prod_{n\geq 1}\left(1-\frac{z^2}{j_{\nu,n}'^2}\right)$$
it follows that
$$2^{\nu}\Gamma(\nu)z^{\frac{3}{2}-\frac{\nu}{2}}J_{\nu}''(\sqrt{z})=z\prod_{n\geq 1}\left(1-\frac{z}{j_{\nu,n}''^2}\right)$$
and
$$-\frac{1}{2}\left(1-\nu+\frac{zJ_{\nu}''(z)}{J_{\nu}'(z)}\right)=\sum_{n\geq1}\frac{z^2}{j_{\nu,n}'^2-z^2}.$$
On the other hand, we know that $\nu\mapsto j_{\nu,n}'$ is increasing on $(0,\infty)$ for each $n\in\mathbb{N}$ fixed (see \cite[p. 236]{nist}), and thus the function
$$\nu\mapsto \sum_{n\geq 1}\frac{1}{j_{\nu,n}'^2-1}=-\frac{1}{2}\left(1-\nu+\frac{J_{\nu}''(1)}{J_{\nu}'(1)}\right)$$ is decreasing on $(0,\infty).$
Consequently, we have that the inequality
$$\sum_{n\geq 1}\frac{1}{j_{\nu,n}'^2-1}\leq 1$$
is valid if and only if $\nu\geq \dot{\nu},$ where $\dot{\nu}$ is the unique root on $(0,\infty)$ of the equation
\begin{equation}\label{eq0}\sum_{n\geq 1}\frac{1}{j_{\nu,n}'^2-1}=1 \ \ \ \Longleftrightarrow\ \ \ (3-\nu)J_{\nu}'(1)+J_{\nu}''(1)=0.\end{equation}
Since $J_{\nu}$ satisfies the Bessel differential equation, it follows that
$$z^2J_{\nu}''(z)+zJ_{\nu}'(z)+(z^2-\nu^2)J_{\nu}(z)=0,$$
and then
$$J_{\nu}''(1)=(\nu^2-1)J_{\nu}(1)-J_{\nu}'(1)=(\nu^2-\nu-1)J_{\nu}(1)+J_{\nu+1}(1),$$
where we used the recurrence relation $zJ_{\nu}'(z)=\nu J_{\nu}(z)-zJ_{\nu+1}(z).$ Consequently, the equation \eqref{eq0} is equivalent to $$(2\nu-1)J_{\nu}(1)+(\nu-2)J_{\nu+1}(1)=0.$$ Now, applying the inequality \cite[Theorem 6.3]{ismail}
$$j_{\nu,1}'^2>\frac{4\nu(\nu+1)}{\nu+2},$$
where $\nu>0,$ it follows that for $n\in\{2,3,\dots\}$ we have $j_{\nu,n}'>{\dots}>j_{\nu,1}'>1$ if $\nu>(-3+\sqrt{41})/8\simeq0.4253{\dots}.$
Thus, by applying Lemma \ref{lem} the assertion of the theorem follows.
\end{proof}

\begin{proof}[\bf Proof of Theorem \ref{thB}]
Let us denote by $j_{\nu,n}''$ the $n$th positive zero of the function $J_{\nu}''.$ By using the infinite product representation \cite[p. 340]{skelton}
$$J_{\nu}''(z)=\frac{\left(\frac{z}{2}\right)^{\nu-2}}{4\Gamma(\nu-1)}\prod_{n\geq 1}\left(1-\frac{z^2}{j_{\nu,n}''^2}\right)$$
it follows that
$$2^{\nu}\Gamma(\nu-1)z^{2-\frac{\nu}{2}}J_{\nu}''(\sqrt{z})=z\prod_{n\geq 1}\left(1-\frac{z}{j_{\nu,n}''^2}\right)$$
and
$$-\frac{1}{2}\left(2-\nu+\frac{zJ_{\nu}'''(z)}{J_{\nu}''(z)}\right)=\sum_{n\geq1}\frac{z^2}{j_{\nu,n}''^2-z^2}.$$
On the other hand, we know that $\nu\mapsto j_{\nu,n}''$ is increasing on $(1,\infty)$ for each $n\in\mathbb{N}$ fixed (see \cite{mercer,wong}), and thus the function
$$\nu\mapsto \sum_{n\geq 1}\frac{1}{j_{\nu,n}''^2-1}=-\frac{1}{2}\left(2-\nu+\frac{J_{\nu}'''(1)}{J_{\nu}''(1)}\right)$$ is decreasing on $(1,\infty).$
Consequently, we have that the inequality
$$\sum_{n\geq 1}\frac{1}{j_{\nu,n}''^2-1}\leq 1$$
is valid if and only if $\nu\geq \nu^{\ast},$ where $\nu^{\ast}$ is the unique root on $(1,\infty)$ of the equation
\begin{equation}\label{eq1}\sum_{n\geq 1}\frac{1}{j_{\nu,n}''^2-1}=1 \ \ \ \Longleftrightarrow\ \ \ (4-\nu)J_{\nu}''(1)+J_{\nu}'''(1)=0.\end{equation}
Since $J_{\nu}$ satisfies the Bessel differential equation, it follows that
$$z^2J_{\nu}'''(z)+3zJ_{\nu}''(z)+(z^2+1-\nu^2)J_{\nu}'(z)+2zJ_{\nu}(z)=0,$$
and then
$$J_{\nu}'''(1)=(1-3\nu^2)J_{\nu}(1)+(\nu^2+1)J_{\nu}'(1)=(\nu^3-3\nu^2+\nu+1)J_{\nu}(1)-(\nu^2+1)J_{\nu+1}(1),$$
where we used the recurrence relation $zJ_{\nu}'(z)=\nu J_{\nu}(z)-zJ_{\nu+1}(z).$ Consequently, the equation \eqref{eq1} is equivalent to $$(2\nu^2-2\nu-3)J_{\nu}(1)=(\nu^2+\nu-3)J_{\nu+1}(1).$$ Now, applying the inequality \cite[Theorem 8.1]{ismail}
$$j_{\nu,1}''^2>\frac{4\nu(\nu-1)}{\nu+2},$$
where $\nu>1,$ it follows that for $n\in\{2,3,\dots\}$ we have $j_{\nu,n}''>{\dots}>j_{\nu,1}''>1$ if $\nu>(5+\sqrt{19})/8\simeq1.5687{\dots}.$
Thus, by applying Lemma \ref{lem} the assertion of the theorem follows.
\end{proof}

\begin{proof}[\bf Proof of Theorem \ref{thC}]
Similarly, as in the proof of the previous theorem, let us denote by $j_{\nu,n}'''$ the $n$th positive zero of the function $J_{\nu}'''.$ By using the infinite product representation \cite[p. 340]{skelton}
$$J_{\nu}'''(z)=\frac{\left(\frac{z}{2}\right)^{\nu-3}}{8\Gamma(\nu-2)}\prod_{n\geq 1}\left(1-\frac{z^2}{j_{\nu,n}'''^2}\right)$$
it follows that
$$2^{\nu}\Gamma(\nu-2)z^{\frac{5}{2}-\frac{\nu}{2}}J_{\nu}'''(\sqrt{z})=z\prod_{n\geq 1}\left(1-\frac{z}{j_{\nu,n}'''^2}\right)$$
and
$$-\frac{1}{2}\left(3-\nu+\frac{zJ_{\nu}''''(z)}{J_{\nu}'''(z)}\right)=\sum_{n\geq1}\frac{z^2}{j_{\nu,n}'''^2-z^2}.$$
On the other hand, we know that $\nu\mapsto j_{\nu,n}'''$ is increasing on $(2,\infty)$ for each $n\in\mathbb{N}$ fixed (see \cite{koko,lorch}), and thus the function
$$\nu\mapsto \sum_{n\geq 1}\frac{1}{j_{\nu,n}'''^2-1}=-\frac{1}{2}\left(3-\nu+\frac{J_{\nu}''''(1)}{J_{\nu}'''(1)}\right)$$ is decreasing on $(2,\infty).$
Consequently, we have that the inequality
$$\sum_{n\geq 1}\frac{1}{j_{\nu,n}'''^2-1}\leq 1$$
is valid if and only if $\nu\geq \nu^{\star},$ where $\nu^{\star}$ is the unique root on $(2,\infty)$ of the equation
\begin{equation}\label{eq2}\sum_{n\geq 1}\frac{1}{j_{\nu,n}'''^2-1}=1 \ \ \ \Longleftrightarrow\ \ \ (5-\nu)J_{\nu}'''(1)+J_{\nu}''''(1)=0.\end{equation}
Since $J_{\nu}$ satisfies the Bessel differential equation, it follows that
$$z^2J_{\nu}''''(z)+5zJ_{\nu}'''(z)+(z^2+4-\nu^2)J_{\nu}''(z)+4zJ_{\nu}'(z)+2J_{\nu}(z)=0,$$
and then
$$J_{\nu}''''(1)=(\nu^4+9\nu^2-2)J_{\nu}(1)-(6\nu^2+4)J_{\nu}'(1)=(\nu^4-6\nu^3+9\nu^2-4\nu-2)J_{\nu}(1)+(6\nu^2+4)J_{\nu+1}(1).$$
Consequently, the equation \eqref{eq2} is equivalent to $$(2\nu^3-7\nu^2+3)J_{\nu}(1)+(\nu^3+\nu^2+\nu-1)J_{\nu+1}(1)=0.$$ Now, taking into account that the function $\nu\mapsto j_{\nu,1}'''$ is increasing on $(2,\infty)$ it follows that for $\nu>3$ we have $j_{\nu,1}'''>j_{3,1}'''\simeq 1.3762{\dots}>1.$ Thus, for $n\in\{2,3,\dots\}$ we have $j_{\nu,n}'''>{\dots}>j_{\nu,1}'''>1$ if $\nu>3,$ and applying again Lemma \ref{lem} the assertion of the theorem follows. We would like to mention here that we approximated the zero $j_{3,1}'''$ by using the mathematical software Matlab by taking into account that $j_{3,1}'''$ is actually the first positive zero of the equation $$\left((1-\nu)z^2+\nu^3-3\nu^2+2\nu\right)J_{\nu}(z)=\left((2+\nu^2)z-z^3\right)J_{\nu+1}(z)$$ when $\nu=3.$
\end{proof}

\begin{proof}[\bf Proof of Theorem \ref{thD}] Let us consider the power series
$$\frac{2^{\nu}\Gamma(\nu+1)}{Q(\nu)z^{\nu}}\left(az^2J_{\nu}''(z)+bzJ_{\nu}'(z)+cJ_{\nu}(z)\right)=
\sum_{n\geq0}\frac{(2n+\nu)(2n+\nu-1)a+(2n+\nu)b+c}{4^{n}n!{(\nu+1)}_n}(-1)^nz^{2n},$$
where $(a)_n=a(a+1)\dots(a+n-1)=\Gamma(a+n)/\Gamma(a).$ By using the fact that for $\tau>0$ the quotient $\log\Gamma(n+\tau)/(n\log n)$ tends to $1$
as $n$ tends to infinity, we obtain that the growth order of the above entire function is the following
$$\rho=\lim_{n\to\infty}\frac{n\log n}{n\log4+\log\Gamma(n+1)+\log\Gamma(n+\nu+1)-\log((2n+\nu)(2n+\nu-1)a+(2n+\nu)b+c)}=\frac{1}{2}.$$
Thus, if $\lambda_{\nu,n}$ denotes the $n$th positive zero of the function $z\mapsto az^2J_{\nu}''(z)+bzJ_{\nu}'(z)+cJ_{\nu}(z),$ then by applying Hadamard's theorem \cite[p. 26]{lev} we obtain
$$az^2J_{\nu}''(z)+bzJ_{\nu}'(z)+cJ_{\nu}(z)=\frac{Q(\nu)z^{\nu}}{2^{\nu}\Gamma(\nu+1)}\prod_{n\geq 1}\left(1-\frac{z^2}{\lambda_{\nu,n}^2}\right),$$
and consequently
$$\frac{2^{\nu}\Gamma(\nu+1)}{Q(\nu)}z^{1-\frac{\nu}{2}}\left(azJ_{\nu}''(\sqrt{z})+b\sqrt{z}J_{\nu}'(\sqrt{z})+cJ_{\nu}(\sqrt{z})\right)=z\prod_{n\geq 1}\left(1-\frac{z}{\lambda_{\nu,n}^2}\right),$$
$$-\frac{1}{2}\left(-\nu+z\cdot\frac{az^2J_{\nu}'''(z)+(2a+b)zJ_{\nu}''(z)+(b+c)J_{\nu}'(z)}{az^2J_{\nu}''(z)+bzJ_{\nu}'(z)+cJ_{\nu}(z)}\right)=
\sum_{n\geq1}\frac{z^2}{\lambda_{\nu,n}^2-z^2}.$$
Here we used the fact that when $\nu\geq \overline{\nu},$ where $\overline{\nu}=\max\{0,\nu_0\}$ and $\nu_0$ is the largest root of the quadratic $Q(\nu)=a\nu(\nu-1)+b\nu+c,$ the zeros of the function $z\mapsto az^2J_{\nu}''(z)+bzJ_{\nu}'(z)+cJ_{\nu}(z)$ are real, according to \cite[Theorem 7.1]{ismail}.

On the other hand, we know cf. \cite[Theorem 1]{mercer} that for $a,b,c\in\mathbb{R}$ such that $c=0$ and $b\neq a$ or $c>0$ and $b>a$ we have that $\nu\mapsto\lambda_{\nu,n}$ is increasing on $(0,\infty)$ for fixed $n\in\mathbb{N}.$ Consequently, under the assumptions, the function
$$\nu\mapsto-\frac{1}{2}\left(-\nu+\frac{aJ_{\nu}'''(1)+(2a+b)J_{\nu}''(1)+(b+c)J_{\nu}'(1)}{aJ_{\nu}''(1)+bJ_{\nu}'(1)+cJ_{\nu}(1)}\right)=
\sum_{n\geq1}\frac{1}{\lambda_{\nu,n}^2-1}$$
is decreasing on $(0,\infty).$ Thus, the inequality
$$\sum_{n\geq 1}\frac{1}{\lambda_{\nu,n}^2-1}\leq 1$$
is valid if and only if $\nu\geq \nu^{\circ},$ where $\nu^{\circ}$ is the unique root on $(\overline{\nu},\infty)$ of the equation
$$\sum_{n\geq 1}\frac{1}{\lambda_{\nu,n}^2-1}=1 \ \ \ \Longleftrightarrow\ \ \ aJ_{\nu}'''(1)+(4a-a\nu+b)J_{\nu}''(1)+(3b+c-b\nu)J_{\nu}'(1)-(\nu-2)cJ_{\nu}(1)=0.$$
By using from the above proofs the expressions for the values $J_{\nu}'''(1),$ $J_{\nu}''(1)$ and $J_{\nu}'(1),$ the above equation is equivalent to
$$(2a\nu^2-2a\nu+2b\nu-3a-b+2c)J_{\nu}(1)=(a\nu^2+a\nu-b\nu-3a+2b+c)J_{\nu+1}(1).$$

Finally, by using the inequalities in \eqref{abc} together with \cite[Eq. (8.2)]{ismail}
$$\lambda_{\nu,1}>\frac{4(\nu+1)Q(\nu)}{Q(\nu)+4a\nu+2a+2b},$$
it follows that for $n\in\{2,3,\dots\}$ we have $\lambda_{\nu,n}>{\dots}>\lambda_{\nu,1}>1,$ and using Lemma \ref{lem} the proof is done.
\end{proof}

\end{document}